\documentclass[12pt]{article} 
\usepackage[sectionbib]{natbib}
\usepackage{array,epsfig,rotating}
\usepackage[]{hyperref}  
\usepackage{sectsty, secdot}
\sectionfont{\fontsize{12}{14pt plus.8pt minus .6pt}\selectfont}
\renewcommand{\theequation}{\thesection\arabic{equation}}
\subsectionfont{\fontsize{12}{14pt plus.8pt minus .6pt}\selectfont}

\usepackage{amsmath}
\usepackage{amssymb}
\usepackage{amsfonts}
\usepackage{multirow}
\usepackage{amsthm}

\usepackage{float}
\usepackage{ulem}
\usepackage{bm}
\usepackage{url}
\usepackage{color} 
\usepackage{lscape}
\usepackage{latexsym}
\usepackage{enumerate}
\usepackage[shortlabels]{enumitem}

\newcommand\bt {\mathbf t}

\newcommand\bx {\mathbf x}
\newcommand\by {\mathbf y}

\newcommand\bW {\mathbf W}

\newcommand\bZ {\mathbf Z}

\newcommand\itA {{\mathcal{A}}}

\newcommand\itE {{\mathcal{E}}}
\newcommand\itF {{\mathcal{F}}}
\newcommand\itG {{\mathcal{G}}}
\newcommand\itH {{\mathcal{H}}}
\newcommand\itI {{\mathcal{I}}}

\newcommand\itL {{\mathcal{L}}}

\newcommand\itN {{\mathcal{N}}}

\newcommand\balfa {\mbox{\boldmath $\alpha$}}

\newcommand\bbe {\mbox{\boldmath $\beta$}}

\newcommand\bkapa {\mbox{\boldmath $\kappa$}}
\newcommand\bkapach {\mbox{\footnotesize\boldmath $\kappa$}}

\newcommand\bmu {\mbox{\boldmath $\mu$}}

\newcommand\btau {\mbox{\boldmath $\tau$}}
\newcommand\btauch {\mbox{\footnotesize\boldmath $\tau$}}

\newcommand\bGa {\mbox{\boldmath $\Gamma$}}
\newcommand\bGach {\mbox{\scriptsize\boldmath $\Gamma$}}

\newcommand\bSi {\mbox{\boldmath $\Sigma$}}

\newcommand\wpsi {\widehat{\psi}}
\newcommand\wphi {\widehat{\phi}}

\newcommand\wlam {\widehat{\lambda}}

\newcommand\wbGa  {\widehat{\bGa}}

\newcommand\witL {\widehat{\itL}}

\newcommand\wtphi {\widetilde{\phi}}
\newcommand\wtpsi {\widetilde{\psi}}
\newcommand\wtlam {\widetilde{\lambda}}

\def\real{\mathbb{R}}
\def\natu{\mathbb{N}}

\newcommand{\HSo}{{\itH_{\smooth}^0}}
\newcommand{\Hdo}{{\itH_{d_n}^0}}

\newcommand{\esp}{\mathbb E}
\newcommand{\prob}{\mathbb{P}}
\newcommand{\cov}{\mbox{\rm Cov}}
\newcommand{\wcov}{\widehat{\cov}}
\newcommand{\wvar}{\widehat{\var}}
\newcommand{\corr}{\mbox{\rm Corr}}

\newcommand{\sd}{\mbox{\textsc{sd}}}
\newcommand{\var}{\mbox{\rm Var}}

\newcommand{\convpp}{ \buildrel{a.s.}\over\longrightarrow}

\newcommand{\trasp}{^\top}

\def\median{\mathop{\mbox{median}}}

\def\argmax{\mathop{\mbox{argmax}}}

\newcommand\noi{\noindent}
\def\dst{\displaystyle}

\def\square{\ifmmode\sqr\else{$\sqr$}\fi}
\def\sqr{\vcenter{
         \hrule height.1mm
         \hbox{\vrule width.1mm height2.2mm\kern2.18mm
\vrule width.1mm}
         \hrule height.1mm}}

\newcommand{\Lr}{\itL_{\rob}}
\newcommand{\Lbtr}{\itL_{\btauch,\rob}}
\newcommand{\Ltr}{\itL_{\tau,\rob}}
\newcommand{\Ltrn}{\itL_{\tau_n,\rob}}

\newcommand{\wLtr}{\witL_{\tau,\rob}}
\newcommand{\wLtrn}{\witL_{\tau_n,\rob}}
\newcommand{\Ga}{{\Gamma}}

\newcommand\gk {\mbox{\footnotesize\sc gk}}

\newcommand{\rob}{{\mbox{{\footnotesize  \sc r}}}}

\newcommand\smooth {\mbox{{\footnotesize\sc s}}}

\def\argmax{\mathop{\mbox{argmax}}}

\newtheorem{theorem}{Theorem}
\newtheorem{lemma}{Lemma}

\newtheorem{proposition}{Proposition}
\theoremstyle{definition}

\newtheorem{remark}{Remark}

\begin{document}


   \title{Robust smoothed canonical correlation analysis  for functional data}
\author{Graciela Boente$^{a,c}$ and Nadia L. Kudraszow$^{b,c}$\\ 
 \small $^a$ Universidad de Buenos Aires   $^b$ Universidad Nacional de La Plata $^c$ CONICET, Argentina }
\date{}

\maketitle


\begin{abstract}
This paper provides robust estimators for the first canonical correlation and directions of random elements on Hilbert separable spaces by using robust association and scale measures  combined with basis expansion and/or penalizations as a regularization tool. Under regularity conditions, the resulting estimators are consistent.

\end{abstract}

\def\thefigure{\arabic{figure}}
\def\thetable{\arabic{table}}

\renewcommand{\theequation}{\thesection.\arabic{equation}}

\fontsize{12}{14pt plus.8pt minus .6pt}\selectfont

\setcounter{section}{0} 
\setcounter{equation}{0} 

\section{Introduction}

In recent years,  data collected in the form of functions or curves received considerable attention in such fields of applications as chemometrics,
economics, environmental studies, image recognition, spectroscopy, and many others. These data are known in the literature as functional data, see \cite{Ramsay} for a complete overview. In general, the observations are considered random elements of some functional space and, in this context, many statistical modelling problems result best described. This gives rise to the extension of some classic concepts of multivariate data analysis, such as dimension reduction techniques and  particularly those based on projections.  

This paper is concerned with canonical correlation analysis, where data consist of pairs of random curves. The aim of this analysis is to identify and quantify the relation between the observed functions. Under a Gaussian model, \cite{Leurgans}  showed that the natural extension of multivariate estimators to the functional scenario fails, which motivates the need to introduce a regularization technique that involves smoothing  through a penalty term. 
 Besides, \cite{He} provided conditions that ensure the existence and proper definition of canonical directions and correlations for processes that support a Karhunen--Lo\'eve expansion, while \cite{Cupidon} derived the asymptotic distribution of correlations and regularized functional canonical variations.  
 An alternative way to get around the  ill--posed problem  related to functional canonical correlation is to use a finite basis expansion. Proposals based on this approach were  discussed in  \cite{He04} and \cite{Ramsay}. More precisely, these authors proposed  to perform   a regularization step projecting the observed curves  on a finite number of basis functions, before computing the smooth canonical correlations and directions in the basis expansion domain.   

All these papers use  the Pearson correlation as   measure of the association between the observed functions. However, it is known the Pearson correlation is  sensitive to atypical observations and this sensitivity is inherited by the procedures based on it (see \cite{Taskinen}). To our knowledge, when considering the analysis of functional canonical correlation, the only proposal of estimators  resistant to anomalous observations is that studied by \cite{Alvarez}, where the regularization was implemented by projecting random processes on a finite number of functions in an orthonormal basis. 

The aim of this paper is to introduce consistent robust estimators of the canonical correlation analysis in the functional data setting but where regularization is based on a roughness penalty.   The paper is organized as follows. In Section \ref{sec:prelim}, we state some notation and preliminary definitions and we briefly describe the two classical approaches for regularized functional canonical correlation analysis. Section \ref{sec:robprop} presents our robust proposals, while their consistency is studied in Section \ref{sec:consis}. Some final comments are given in Section \ref{sec:conclusion}. All
proofs are deferred to the Appendix.

\section{Preliminaries}{\label{sec:prelim}}
Let $\itH$ be  a separable Hilbert space with inner product $\left\langle \cdot,\cdot\right\rangle$ and norm  $\left\| \cdot \right\|^2= \left\langle \cdot,\cdot\right\rangle$. Let $(X,Y)$  be a random element of the Hilbert space $\itH\times \itH$ defined in a probability space $(\Omega,{\itA}, \prob)$. In the product space $\itH\times\itH$, we define  the usual inner product  ${\langle (u_1,v_1),(u_2,v_2)\rangle}_{\itH\times \itH}={\langle u_1,u_2\rangle}+{\langle v_1,v_2\rangle}$. 
When $(X,Y)\trasp$ has finite second moment, i.e., $\esp(\|X\|^2+\|Y\|^2)<\infty$,
we denote as  $\bGa_{XX}:\itH\to\itH$, $\bGa_{YY}:\itH\to\itH$, $\bGa_{XY}:\itH\to\itH$ and $\bGa_{YX}:\itH\to\itH$,  the covariance and cross--covariance operators, respectively. More precisely, for any $u_1, u_2 \in \itH$, $v \in \itH$, we have that $\cov(\langle u_1,X\rangle,\langle u_2,X\rangle)=\langle u_1,\bGa_{XX} u_2 \rangle$, $\cov(\langle u_1,X\rangle,\langle v,Y\rangle)=\langle u_1,\bGa_{XY} v \rangle$ and  similarly for  $ \bGa_{YY}  $ and $  \bGa_{YX} $.

\subsection{The classical approaches}\label{classical}
Canonical correlation analysis,
which was originally developed for multivariate data, has been successfully extended to accommodate functional data  by \cite{Leurgans} as follows. 

Assume that the observed data  $\{(X_i , Y_i )\trasp, i = 1, \dots, n \}$ are independent realizations of a bivariate stochastic process $(X,Y)\trasp\in \itH\times\itH$.   When $(X,Y)\trasp$ has finite second moments, a non-smooth approach to the problem of functional canonical correlation is to search for functions $u$ and $v$ in $\itH$ such that the linear combinations $\langle u,X\rangle$ and $\langle v,Y\rangle$ have maximum squared correlation, that is, the objective is to find $u\neq 0$, $v\neq 0$ that maximize
\begin{equation}
\label{eq:ITL}
\itL(u,v)=\corr^2(\left\langle u,X\right\rangle,\left\langle v,Y\right\rangle)= \frac{\left\langle u, \bGa_{XY}v\right\rangle^2}{\left\langle u, \bGa_{XX}u\right\rangle\left\langle v, \bGa_{YY}v\right\rangle}\;,
\end{equation}
where  the ratio  ${\langle u,\bGa_{XY}v\rangle^2}/ \left(\langle u,\bGa_{XX}u\rangle\;\langle v,\bGa_{YY}v\rangle\right)  $   equals $0$ when $\langle u,\bGa_{XX}u\rangle=0$ or $\langle v,\bGa_{YY}v\rangle=0$.
In particular, \cite{Leurgans} considered the case $\itH = L^2(\itI)$ and assumed that there are two bases of $\itH$ composed of the functional canonical coordinates,  which are a generalization of the vector canonical coordinates,  that  ensure  the existence of a solution to the non-smooth approach.

\cite{Leurgans} proved that it is not possible to consider a sample version of the problem of maximizing $\itL(u,v)$. Therefore, they proposed to estimate the first canonical variables by maximizing, restricted to not null `smooth elements' of $\itH$, the estimated canonical correlation penalized by a `penalty operator'. 

As mentioned in the Introduction,  two possibilities may be considered to introduce regularization. 
One approach is to consider, as in \cite{Leurgans},  a roughness penalty which gives a measure of the smoothness of a function.  The other point of view, considers  a sieves approximation eventually combined with a penalty term.  We will briefly review both methods.

Let $D:\itH_{\smooth}\rightarrow\itH$ be a linear operator, which we will refer to as the \textsl{differentiator}, here $\itH_{\smooth}$ is the subset of \textsl{smooth elements} of $\itH$, i.e.,  $u\in\itH_{\smooth}$ if $\|Du\| <\infty$. 
Using $D$, we define the symmetric positive semi-definite bilinear form $\left\lceil \cdot,\cdot\right\rceil :\itH_{\smooth}\times \itH_{\smooth}\rightarrow\real$,
where $\left\lceil u, v\right\rceil=\left\langle D u,D v\right\rangle$. The \textsl{penalization operator} is then defined as
$\Psi:\itH_{\smooth}\rightarrow\real$,	$\Psi(u) = \left\lceil u, u\right\rceil$, and the penalized inner product as $\langle u, v\rangle_{\tau} =\langle u,v\rangle+\tau\left\lceil u, v\right\rceil$.

\begin{remark}\label{remark:remark1}  The most common setting for functional data corresponds to the situation where $\itH=L^2(\itI)$ and 
$$\itH_{\smooth} =\{ u \in L^2(\itI),\; u \text{ is }\newline \text{twice differentiable, and }
\int_{\itI}{u^{\prime\,\prime}(t)}^2 dt <\infty\}\,.$$ 
In this case, it is usual to consider $D u = u^{\prime\,\prime}$ and $\left\lceil u,v\right\rceil =
\int_{\itI}u^{\prime\,\prime}(t)\, v^{\prime\,\prime}(t) dt$, so that $\Psi(u) =
\int_{\itI}{u^{\prime\,\prime}(t)}^2dt$. \hfill $\clubsuit$
\end{remark}

\vskip0.2in

Denote as $\itH_{\smooth}^0:=\{u \in\itH_{\smooth}:\,u\neq 0\}$. Given $u$ and $v$ in $\itH_{\smooth}^0$,   \cite{Leurgans} defined  the population penalized squared correlation, $\itL_{\btau}(u,v)$,   as 
\begin{eqnarray*}
\itL_{\btau}(u,v)&=&\frac{\cov^2(\left\langle u,X\right\rangle,\left\langle v,Y\right\rangle)}{\left\{\var(\left\langle u,X\right\rangle)+\tau_1\Psi(u)\right\}\left\{\var(\left\langle v,Y\right\rangle)+\tau_2\Psi(v)\right\}}\\
&=&\frac{\left\langle u, \bGa_{XY}v\right\rangle^2}{\left\{\left\langle u, \bGa_{XX}u\right\rangle+\tau_1\Psi(u)\right\}\left\{\left\langle v, \bGa_{YY}v\right\rangle+\tau_2\Psi(v)\right\}}\;,
\end{eqnarray*}
where $\btau=(\tau_1,\tau_2)$. 
 The so--called \textsl{Smoothed Canonical Correlation Analysis} (SCCA) by \cite{Leurgans}, correspond to   maximizing $\itL_{\btau}(u,v)$ over $u,v\in  \itH_{\smooth}^0$. 
 In this way, for the sample $\{(X_i, Y_i)\trasp, i = 1, \dots, n \}$  these authors proposed to carry out  SCCA  by  replacing the population quantities by their sample counterparts, that is, by maximizing the penalized squared sample correlation
\begin{eqnarray}
\witL_{\btau}(u,v)&=&\frac{\wcov^2(\left\langle u,X\right\rangle,\left\langle v,Y\right\rangle)}{\left(\wvar(\left\langle u,X\right\rangle)+\tau_1\Psi(u)\right)\left(\wvar(\left\langle v,Y\right\rangle)+\tau_2\Psi(v)\right)}
\label{eq:clwitl}\\
&=&\frac{\left\langle u, \wbGa_{XY}v\right\rangle^2}{\left(\left\langle u, \wbGa_{XX}u\right\rangle+\tau_1\Psi(u)\right)\left(\left\langle v, \wbGa_{YY}v\right\rangle+\tau_2\Psi(v)\right)}\;,
\nonumber
\end{eqnarray}
where $\wcov$ and $\wvar$ stand for the sample covariance and variance, computed replacing the corresponding bivariate or univariate distributions with the empirical one, respectively, while $ \wbGa_{XX}$, $ \wbGa_{YY}$ and $\wbGa_{XY}$ stand for the sample covariance and cross--covariance operators, respectively.

 As mentioned in the Introduction, to address the dimensionality problems of functional canonical correlation,  \cite{He04} and \cite{Ramsay} propose an alternative to SCCA by means of dimension reduction techniques, that is, following a sieves approach. More precisely, these authors implement regularization by first projecting the sample's curves on a finite number of elements of an orthonormal basis.  In this way, given $\{\xi_i\}_{i\geq 1}$  a suitable orthonormal basis  for $\itH$, we will denote by $\itH_d$ the subspace of $\itH$ spanned by $\{\xi_1,\dots,\xi_d\}$. Then, if we take $d=d_n$ such that $d_n\to \infty$, the sequence of increasing subspaces $\itH_{d_n}$ approximates  $\itH$. From now on, we assume the basis elements are smooth and so $\itH_{d_n}^0:=\{u \in\itH_{d_n}:\,u\neq 0\}\subset\itH_{\smooth}^0$. For simplicity, we only consider the case where the same basis is used to approximate the  both canonical direction estimators.

For the sample $\{(X_i, Y_i)\trasp, i = 1, \dots, n \}$,  \cite{Ramsay} defined the SCCA restricted to the basis expansion domain as the maximization of $\witL_{\btau}(u,v)$ over  $\itH_{d}^0\times\itH_{d}^0$. Let $\balfa=(\alpha_1,\ldots,\alpha_{d})\trasp$ and $\bbe=(\beta_1,\ldots,\beta_{d})\trasp$ be the coefficients' vectors of $u$ and $v$ in the considered basis and denote as $\bx=\left(\langle X,\xi_1\rangle, \ldots, \langle X,\xi_{d}\rangle\right)\trasp$ and $\by=(\langle Y,\xi_1\rangle, \ldots, \langle Y,\xi_{d}\rangle)\trasp$. It is easily seen that, in the basis
expansion domain, the SCCA of the given data is {carried} out by maximizing, over $\balfa,\bbe{\neq \mathbf{0}}$, 
\begin{equation}\label{sievescl}
\witL_{\btau}^d (\balfa,\bbe)=\frac{\wcov^2(\balfa\trasp \bx,\bbe\trasp \by)}{\left(\wvar(\balfa\trasp \bx)+\tau_1\dst\sum_{i,j}^d\alpha_i\, \alpha_j\, \left\lceil \xi_i,\xi_j\right\rceil\right)\left(\wvar(\bbe\trasp \by)+\tau_2\dst \sum_{i,j}^d\beta_i\, \beta_j\,\left\lceil \xi_i,\xi_j\right\rceil\right)}.
\end{equation}
The maximizers of \eqref{sievescl} are the coefficient's vectors of the estimated leading canonical directions in the considered basis.  

Some of the most frequently used bases for functional data are the Fourier, polynomial, splines and wavelet bases. It could also be taken a data-driven basis such as the one composed of the eigenfunctions of the covariance operators. The number of basis elements, $d$, should be chosen large enough to ensure that the regularization is controlled by the choice of the smoothing parameter $\btauch$ rather than that of dimensionality $d$.

\subsection{Co--association measures}\label{measures}
As it is well known,  the estimators obtained maximizing  $\witL_{\btau}(u,v)$  are very sensitive to the presence of outliers, since they are based on the sample version of the covariance operators. This suggests that more resistant association measures are needed to get reliable estimations, see for instance, \cite{Alfons} who introduces  robust canonical correlation estimators for multivariate data and provides a deep discussion on bivariate association measures. Association measures are an alternative to and include the Pearson correlation. In our setting, we seek for robust alternatives to the covariance between two random variables since we are penalizing   the two variances appearing in the denominator of \eqref{eq:clwitl}. Clearly, a resistant measure can be constructed from a robust association measure and a robust scale estimator. However, other possible choices can be considered. We first give a   definition, that provides, a general  framework to  robust counterparts of the usual covariance.

Given two univariate random variables  $U$ and $V$, let $F_{(U,V)}$, $F_U$ and $F_V$ stand for the distributions of $(U,V)\trasp$, $U$ and $V$ respectively.  A \textsl{bivariate  co--association measure}  $\gamma$ between $U$ and $V$, denoted $\gamma(F_{(U,V)})$, is a functional defined over the space of bivariate distributions  such that 
\begin{description}
	\item[(i)]  $\gamma(F_{(U,V)})=\gamma(F_{(V,U)})$,
	\item[(ii)]  $\gamma(F_{(aU+b,cV+d)})= a\,c\,\gamma(F_{(U,V)})$, where $a, b, c$ and $d$ are real constants.
\end{description}
To avoid for burden notation,   we   write  $\gamma(U,V)$ instead of $\gamma(F_{(U,V)})$, from now on.

Furthermore, if a bivariate  co--association measure  $\gamma$ also  satisfies the condition 
\begin{description}
	\item[(iii)] $\gamma^2(U,V)\le \gamma(U,U)\;\gamma(V,V)$,
\end{description}
a   measure of association may be defined as $\rho\left( U,V\right)= {\gamma \left(U,V \right)}/{\sqrt{\gamma \left( U,U\right)\gamma\left( V,V\right)}}$.
 Clearly, the covariance between two random variables is a co--association measure that satisfies \textbf{(i)}-\textbf{(iii)} and its related association measure is the Pearson correlation.

As mentioned above, to provide a robust counterpart of \eqref{eq:clwitl}, robust scale estimators are also needed. To recall  the definition of a scale functional, denote  $\itG$ the set of all univariate distributions. A scale functional $\sigma:\itG\to [0,+\infty)$ is a   location invariant and scale equivariant functional, that is,   $\sigma(F_{aU+b})=|a|\sigma(F_U)$,
for all real numbers $a$ and $b$ (see \cite{Yohai}).  Two well known examples of scale functionals are the standard deviation
and the median absolute deviation about the  median, $\textsc{mad}(F_U)= c \,\median\left(|U-\median(U)| \right)$. The normalization constant $c$, used in the \textsc{mad}, can be chosen so that its
empirical or sample version is consistent for a scale parameter of interest.  Typically, one chooses  $c=1/\Phi^{-1}(0.75)$ so that the \textsc{mad} equals the standard deviation at a normal distribution. More generally, any $M-$scale estimator can be calibrated to provide   Fisher--consistent estimators at the normal distribution, that is, $\sigma(\Phi)=1$, with $\Phi$ the standard normal distribution. As above, when there is no confusion, we will denote $\sigma(U)$ instead of $\sigma\left(F_U\right)$.

Given a  bivariate co--association functional  $\gamma$ and a  scale functional $\sigma$, one can define  the related association measure  $\rho$ as 
$
\rho\left( U,V\right)=  \gamma \left( U,V \right)/\{\sigma \left( U\right)\sigma\left( V\right)\}
$, 
if $\rho^2\left( U,V \right)\le 1$, for any two univariate variables $U$ and $V$.
Conversely, given an association measure $\rho$ and  a  scale functional $\sigma$, the related co--association is given by
$\gamma \left( U,V \right)= \rho\left( U,V \right)\;\sigma \left( U\right)\sigma\left( V\right)$.

Examples of such association measures can be constructed from   a bivariate robust scatter functional  $\bW=\bW(U,V)$, which provides a more resistant alternative to the classical covariance matrix $\bSi=\cov(U,V)$. The association measure induced by a bivariate scatter matrix $\bW$ is given by
\begin{eqnarray}\label{eq:rhocov}
\rho(U,V)= \frac{\bW_{12}(U,V)}{\left\{\bW_{11}(U,V)\bW_{22}(U,V)\right\}^{\frac 12}}, 
\end{eqnarray}
where $\bW_{ij}(U,V)$ is the $(i,j)-$th element of the scatter matrix $\bW(U,V)$.
One possible choice for  $\bW(U,V)$ is the $M-$scatter estimator  defined by  \cite{Maronna}, since it provides an efficient estimator which is also highly robust in the bivariate case. Another possible choice is to consider the orthogonalized Gnanadesikan--Kettenring covariance  proposed by \cite{MaronnaZ}.  When using $M$-estimators or the orthogonalized Gnanadesikan--Kettenring covariance,  the corresponding    co--association measure is defined taking $\gamma(U,V)=\bW_{12}(U,V)$ and the related scale estimators as $\sigma(U)=\sqrt{\bW_{11}(U,V)}$ and  $\sigma(V)=\sqrt{\bW_{22}(U,V)}$. Note that $\rho^2\left( U,V\right)\le 1$, when $\bW$ is positive semi-definite, which is satisfied by both estimators mentioned above.

Taking into account that $\cov(U,V)=({\alpha\beta}/{4})({\sd}^2(U/\alpha+V/\beta)-{\sd}^2(U/\alpha-V/\beta))$ for all $\alpha\neq 0$ and $\beta \neq 0$,
where  $\sd(\cdot)$ stands for  the standard deviation, \cite{Gnanadesikan} define  a family of co--association functionals replacing the standard deviation by a robust scale $\sigma$ and taking $\alpha=\sigma(U)$ and $\beta=\sigma(V)$. More precisely, given a scale functional $\sigma$,  the co--association measure $\gamma^{\star}$ is defined as $ 
\gamma^{\star}(U,V)= {\sigma(U)\,\sigma(V)} \, (\sigma^2_{+}- \sigma^2_{-})/4$ with
\begin{equation}
\sigma^2_{+}=\sigma^2\left(\frac{U}{\sigma(U)}+\frac{V}{\sigma(V)}\right) \qquad \sigma^2_{-}=\sigma^2\left(\frac{U}{\sigma(U)}-\frac{V}{\sigma(V)}\right) \label{eq:sigma+-}
\end{equation}
In order to obtain a highly robust estimator of the correlation between two  real random  variables, the association measure $\rho^{\star}(U,V)$ is defined as in  $\rho^{\star}(U,V)=(\sigma^2_{+}- \sigma^2_{-})/4$. However,  the resulting measure will  not  bounded  between  $-1$  and  $1 $, since the co--association measure does not satisfy \textbf{(iii)}. 
To ensure an association measure in the valid range, \cite{Gnanadesikan} define the association measure $\rho_{\gk}$ as $\rho_{\gk}(U,V)=(\sigma^2_{+}-\sigma^2_{-})/(\sigma^2_{+}+\sigma^2_{-}) $, with $\sigma^2_{+}$ and $\sigma^2_{-}$ defined in \eqref{eq:sigma+-}.
which  lies in
the range $[-1, 1]$, and the related a co--association measure  through $\gamma_{\gk}(U,V)= \sigma(U)\,\sigma(V)\,\rho_{\gk}(U,V)$.

\begin{remark}\label{remark:remark2} We say that  $(U,V)\sim \itE_2(\bmu, \bSi,\varphi)$ if $\bZ=(U,V)\trasp$ is elliptically distributed with location $\bmu$, scatter matrix $\bSi$ and  characteristic generator function $\varphi$, i.e., the characteristic function of   $\bZ$ equals  $\psi_{\bZ}(\bt)=\exp(i \bmu\trasp \bt) \varphi(\bt\trasp \bSi \bt)$. 
As mentioned in Section 2.1 in \cite{Alvarez},   if the robust scatter functional $\bW$  is affine--equivariant,  the  association measure defined in \eqref{eq:rhocov} is Fisher--consistent for elliptical families,  that is,   $\rho(U,V)=\Sigma_{12}/\sqrt{\Sigma_{11}\; \Sigma_{22}}$. 
In particular, the association measure induced by the $M-$scatter estimator  defined by \cite{Maronna} is Fisher--consistent at any elliptical distribution. Furthermore, even when the scatter matrix defined in \cite{MaronnaZ} is not affine equivariant,  the association measure $\rho$ given in \eqref{eq:rhocov} is also  Fisher--consistent at any elliptical distribution.

When the scale function $\sigma(\cdot)$ is calibrated so as to be Fisher--consistent at the normal distribution,  $\gamma^{\star}$ and $\gamma_{\gk} $ are Fisher-consistent at the bivariate normal distribution. When considering elliptical distributed random vectors $(U,V)\sim \itE_2(\bmu, \bSi,\varphi)$, it is well known that for any robust scale functional there exists a constant $c>0$ such that for any $a,b\in \real$,  $\sigma^2(aU+bV)=c\left(a^2\Sigma_{11}+b^2\Sigma_{22}+2\,a\,b \Sigma_{12}\right)$ (see, for instance, \cite{Yohai}). Straightforward arguments allow to show that, in such situation,   $\sigma^2_{+}= 2(1+\Sigma_{12}/\sqrt{\Sigma_{11}\; \Sigma_{22}}) $ and $\sigma^2_{-}= 2(1-\Sigma_{12}/\sqrt{\Sigma_{11}\; \Sigma_{22}})$, so $\rho^{\star}$ and $\rho_{\gk}$ are also Fisher--consistent at  elliptical distributions.  \hfill $\clubsuit$
\end{remark}

\section{Robust  approaches  for smoothed canonical correlation analysis}{\label{sec:robprop}}
Throughout this paper,  we will denote as $P_Z[u]$  the distribution of $\left\langle u,Z\right\rangle$ when $Z\sim P_Z$ and as $P_{(X,Y)}[u,v]$  the joint distribution of $(\left\langle u,X\right\rangle,\left\langle v,Y\right\rangle)\trasp$ when $(X,Y)\trasp\sim P_{(X,Y)}$. Furthermore, given a sample $Z_1,\dots, Z_n$, 
$ P_{n,Z}[u]$ stands for  the empirical distribution of $\left\langle u,Z_{1}\right\rangle, \dots,\left\langle u,Z_n\right\rangle$, while  
$P_{n,(X,Y)}[u,v]$ is the the empirical distribution of the bivariate sample $\left(\left\langle u,X_{i}\right\rangle,\left\langle v,Y_{i}\right\rangle\right)\trasp$, $1\le i\le n$. 

Let $\gamma_\rob$ and $\sigma_\rob$ be robust co--association and scale functionals, respectively, defining  a measure of association, that is,  $ \gamma_\rob^2 \left( U,V \right)\le \sigma_\rob^2 \left( U\right)\sigma_\rob^2\left( V\right)$.  From now on,   $\gamma_{XY}(u,v)=\gamma_{\rob}\left(P_{(X,Y)}[u,v]\right)$, $\sigma_{Z}(u)=\sigma_{\rob}\left(P_Z[u]\right)$ while their sample versions will be denoted as $g_n(u,v) = \gamma_\rob(P_{n,(X,Y)}[u,v])$ and  $s^2_{n,Z}(u) =\sigma^2_\rob(P_{n,Z}[u])$, respectively.  When $\gamma_\rob \left( U,V \right)= \rho_\rob\left( U,V \right)\;\sigma_\rob \left( U\right)\sigma_\rob\left( V\right)$ for some  association measure $\rho_\rob$, we will denote as $r_n(u,v) = \rho_\rob(P_{n,(X,Y)}[u,v])$ and $\rho_{XY}(u,v)=\rho_{\rob}\left(P_{(X,Y)}[u,v]\right)$.  Furthermore, given any $u, v \in \itH$, denote as  $\itL_{\rob}(u,v)$  the robust population squared measure of association between $\left\langle u,X\right\rangle$ and $\left\langle v,Y\right\rangle$ and by $\Lbtr(u,v)$ its smoothed version,  that is, 
$$ \itL_{\rob}(u,v)=\frac{\gamma^2_{XY}(u,v)}{\sigma^2_X(u)\sigma^2_Y(v)} \hspace{.4cm} \mbox{and} \hspace{.4cm}\Lbtr(u,v)=\frac{\gamma_{XY}^2(u,v)}{\{\sigma^2_{X}(u)+\tau_1\Psi(u)\}\{\sigma^2_{Y}(v)+\tau_2\Psi(v)\}} \,,
$$
where we define $\itL_{\rob}(u,v)=0$ when $\sigma^2_X(u)=0$ or $\sigma^2_Y(v)=0$. Note that $ \itL_{\rob}$ is the robust counterpart of $\itL(u,v)$ in \eqref{eq:ITL}.    Moreover, if  $\gamma_\rob$ is related to an association measure $\rho_\rob$  and the scale functional $\sigma_\rob$ as $\gamma_\rob \left( U,V \right)= \rho_\rob\left( U,V \right)\;\sigma_\rob \left( U\right)\sigma_\rob\left( V\right)$, then $ \itL_{\rob}(u,v)= \rho_{\rob}^2\left( P_{(X,Y)}[u,v]  \right)$.  We will  refer to the supremum of $\itL_{\rob}(u,v)$    as   the first or maximum canonical association.

As mentioned in Section \ref{classical}, when the  co--association measure  $\gamma_\rob$  and the scale functional $\sigma_\rob$ are taken as the covariance and the standard deviation, functional canonical correlation is an ill--posed problem and some regularization is needed. Similarly, when considering general co--association and scale functionals, it is not possible to consider a sample version of the problem of maximizing $\itL_{\rob}(u,v)$. More precisely, Proposition 3.1 of  \cite{Alvarez} entails that, when $\mbox{dim}(\itH)=\infty$, there are directions such that the empirical association measure $ \witL_{\rob}(u,v)=  g_n^2(u,v)/\{s_{n,X}^2(u)\,s^2_{n,Y}(v)\}$ equals  one. For that reason,  the proposal given in \cite{Leurgans} can easily be adapted, using  the sample version of $\Lbtr$, to get more stable estimators. To simplify our notation, in what
follows  we avoid the subscript \textsc{r} when defining the canonical directions   functionals and their estimators. 
The  robust canonical direction functionals and their smoothed versions are defined, respectively, as $(\phi_{1},\psi_{1})=\argmax_{u,v\,\in\HSo}\itL_{\rob}(u,v)$ and 
$(\phi_{\btauch,1},\psi_{\btauch,1})=\argmax_{u,v\,\in\HSo}{\Lbtr(u,v)}$.
The sample counterparts of $(\phi_{\btauch,1},\psi_{\btauch,1})$ are obtained using the sample versions of the robust co--association and scale functionals, that is, the smoothed robust canonical correlation estimators are given by
\begin{equation}
(\wphi_{\btauch,1},\wpsi_{\btauch,1})=\argmax_{u,v\,\in\HSo}\frac{g_n^2(u,v)}{\{s_{n,X}^2(u)+\tau_1\Psi(u)\}\{s^2_{n,Y}(v)+\tau_2\Psi(v)\}}=\argmax_{u,v\,\in\HSo}\witL_{\btauch,\rob}(u,v)\,.
\label{eq:FCCest}
\end{equation}

In the same way, the proposal of \cite{Ramsay} based on re\-gu\-larization by means of both orthonormal bases  and a penalization parameter can be easily adapted to be robust maximizing $\witL_{\btauch,\rob}$ over $\itH_d^0$. So, the smoothed robust canonical correlation estimators in the basis expansion domain are given by
\begin{equation}
(\wtphi_{\bkapach,1},\wtpsi_{\bkapach,1})=\argmax_{u,v\,\in\itH_d^0}\witL_{\btauch,\rob}(u,v)\,,
\label{eq:FCCestsieves}
\end{equation}
where $\bkapa=(\btau,d)$.

Note that the above maximizations have no unique solution, any scalar multiplication of a solution is also a solution. For that reason, conditions over the norms of the directions or the variances of the projections are usually imposed in order to achieve identifiability up to a sign.  With this equivalence in mind, we have that   $(\phi_{1},\psi_{1})$ is the pair of leading robust canonical directions of the model, while $(\wphi_{\tau,1},\wpsi_{\tau,1})$ and $(\wtphi_{\bkapach,1},\wtpsi_{\bkapach,1})$, given in \eqref{eq:FCCest} and \eqref{eq:FCCestsieves} respectively,  are  its  estimators. It is worth noticing that an unsmoothed robust version of \eqref{eq:FCCestsieves}, i.e., when $\tau_1=\tau_2=0$ was studied in \cite{Alvarez}.

\section{Consistency}{\label{sec:consis}}
As in \cite{Leurgans}, to derive consistency results for {smoothed robust canonical correlation estimators}, it is enough to consider the special case where $\tau_1 = \tau_2 = \tau$, and from now on we shall confine attention to this case.  We will denote as   $\itN$  the null space of $\left\lceil \cdot ,\cdot\right\rceil$ and as $\itN^{\perp}$   its orthogonal complement. 
The following   assumptions are needed to obtain the desired convergence results.

\begin{enumerate}[label = \textbf{C\arabic*}]
	
	\item \label{A1} There exists a constant $c>0$ and a self--adjoint, positive, compact operator  $\bGa:\itH\times \itH\rightarrow\itH\times \itH$ such that 
	\begin{equation}\label{eq:Gamma}
	\bGa=\left(\begin{array}{cc}
	\bGa_{11}&\bGa_{12}\\
	\bGa_{21}&\bGa_{22}
	\end{array}\right)\,,
	\end{equation} 
	and  for any $u,v\in\itH$, $\sigma^2_X(u)=c\left\langle u,\bGa_{11} u\right\rangle$, $\sigma^2_Y(v)=c\left\langle v,\bGa_{22} v\right\rangle$ and $\gamma_{XY}(u,v)=c\left\langle u,\bGa_{12} v\right\rangle$. Besides,  the eigenfunctions of $\bGa_{11}$ and $\bGa_{22}$ fall in $\itH_{\smooth}$.

	\item\label{A2} There exist functions $\phi_{1}$ and $\psi_{1}$ in $\itH_{\smooth}$  such that, for any $u, v \in \itH$, we have $\itL_{\rob}(u,v)\leq  \itL_{\rob}(\phi_{1},\psi_{1})=\rho_0^2$. Furthermore, there exists  $0 \leq \rho_1 < \rho_0$, such that $\itL_{\rob}(u, v) \leq \rho_1$, for any $u\in \itH$ and $v\in \itH$ such that $\itL_{\rob}(u,\phi_{1}) = \itL_{\rob}(v,\psi_{1}) = 0$. {Furthermore, assume that $\|\phi_{1}\|=1$ and $\|\psi_{1}\|=1$ and that $(\phi_{1},\psi_{1})$ is unique up to change of sign.}
	
	\item\label{A3}  
	\begin{enumerate} 	\item For any $u\in\itN$, $u\neq0$, $\sigma_X(u)\ne 0$ and $\sigma_Y(u)\ne 0$.
		\item	$\itN$ is finite dimensional and  there exists $d > 0$ such that $\Psi(u)=\left\lceil u ,u\right\rceil > d\|u\|^2$ for all $u\in\itN^{\perp}$. 
	\end{enumerate}
\end{enumerate}

Note that in \ref{A1} we may assume without loss of generality that $c=1$ redefining $\bGa$ as $c\; \bGa$. From now on, we will denote 
$\|u\|^2_{1,\tau} =\sigma^2_X(u)+\tau\left\lceil u, u\right\rceil =\langle u,\bGa_{11}u\rangle+\tau\left\lceil u, u\right\rceil  $ and $\| v \|^2_{2,\tau} =\sigma^2_Y(v)+\tau\left\lceil v, v\right\rceil =\langle v,\bGa_{22}v\rangle+\tau\left\lceil v, v\right\rceil  $. Furthermore, let  
\begin{align*}
C_{n,X}&=\sup_{\|{u}\|_{1,\tau{_n}}=1}\left|s_{n,X}^2(u)-\sigma^2_X(u)\right|\,, \\ 
C_{n,Y} &=\sup_{\|{v}\|_{2,\tau{_n}}=1}\left|s_{n,Y}^2(v)-\sigma^2_Y(v)\right|\, ,  \\
C_{n,XY} & =\sup_{\|{u}\|_{1,\tau{_n}}=\|{v}\|_{2,\tau{_n}}=1}\left|g_{n}(u,v)-\gamma_{XY}(u,v)\right|\,.
\end{align*} 

We will also need the following assumption which is related to the convergence of the scale and co--association estimators.
\begin{enumerate}[label = \textbf{C\arabic*}]
	\setcounter{enumi}{3}	
	\item\label{A4}  {The smoothing parameter $\tau=\tau_n\geq 0$ is such that   $\tau_n\rightarrow 0$, $ \max(C_{n,X},C_{n,Y})\convpp 0$ and one of the following hold
		\begin{enumerate}
			\item  $ C_{n,XY}\convpp 0$ as $n\rightarrow \infty$ and there exists a constant $A>0$ such that for any $u,v\in \HSo$, we have that 
			$$  {g_{n}^2(u,v)}/\left(\{s_{n,X}^2(u)+\tau\Psi(u)\}\{s_{n,Y}^2(v)+\tau\Psi(v)\}\right) \le A.$$
			\item The co--association measure is such that $g_{n}(u,v)=r_n(u,v) s_{n,X} (u)s_{n,Y}(v)$ and we have that $\theta_n=\sup_{\|{u}\|=\|{v}\|=1}\left|r_{n}(u,v)-\rho_{XY}(u,v)\right|\convpp 0$. 
	\end{enumerate}}
\end{enumerate}
It is worth noticing that the condition $$  {g_{n}^2(u,v)}/[\{s_{n,X}^2(u)+\tau\Psi(u)\}\{s_{n,Y}^2(v)+\tau\Psi(v)\}]\le A$$ for any $u,v\in \HSo$, clearly holds with $A=1$,  when $\gamma_\rob \left( U,V \right)=\rho_\rob\left(U,V \right)\sigma_\rob \left( U\right)\sigma_\rob\left( V\right)$ for some  association measure $\rho_\rob$, as is the case in the classical setting, since $r_n(u,v)\le 1$.  Note also that  $r_n(u,v)$ y $\rho_{XY}(u,v)$ are scale invariant, so we also have that $\theta_n=\sup_{\|{u}\|_{1,\tau{_n}}=\|{v}\|_{2,\tau{_n}}=1}\left|r_{n}(u,v)-\rho_{XY}(u,v)\right|$.

Assume that \ref{A2} holds and define $$\lambda_0 = \rho_0^2= \sup_{u,v\neq 0}\itL_{\rob}(u, v),$$ 
$$
\lambda_{\tau} = \sup_{u,v\,\in\HSo} \Ltr(u,v) \qquad \text{ and } \qquad\wlam_{\tau}  = \sup_{u,v\,\in\HSo}\wLtr(u, v).$$ As in Section \ref{sec:robprop}, we denote  the maximizing values in each case by  $(\phi_{1},\psi_{1})$, $(\phi_{\tau,1},\psi_{\tau,1})$ and $(\wphi_{\tau,1},\wpsi_{\tau,1})$, respectively. 

The notion of convergence of the  canonical directions estimators, $(\wphi_{\tau,1},\wpsi_{\tau,1})$, to the first population canonical ones, $(\phi_{1},\psi_{1})$,  will be the convergence with respect to the   association measure induced by $\gamma_{\rob}$ and $\sigma_{\rob}$, that is analogous to the $\Ga-$norm convergence defined in \cite{Leurgans}. This convergence means that the canonical variates obtained from $(u_n,v_n)$ for a given random element $(X,Y)\trasp$ behave as those obtained from $(u,v)$ which is a desirable property to hold for the estimated canonical directions. To clarify the convergence to be considered, given $u_1, u_2, v_1, v_2 \in \itH$, define the quantities
$$
\Lr^X(u_1, u_2) =\frac{\gamma^2_\rob(P_{(X,X)}[u_1,u_2])}{\sigma^2_X(u_1)\sigma^2_X(u_2)}\qquad \mbox{and}\qquad
\Lr^Y(v_1, v_2) =\frac{\gamma^2_\rob(P_{(Y,Y)}[v_1,v_2])}{\sigma^2_Y(v_1)\sigma^2_Y(v_2)}\,.
$$
For any pair of sequences $(u_n)_{n\in\natu}\subseteq \itH$, $(v_n)_{n\in\natu}\subseteq \itH$, we say that $(u_n,v_n)$ converges to $(u,v)\in \itH$ in the $\itL_{\rob}-$norm   if $\Lr^X(u, u_n)\to 1$ and $\Lr^Y(v, v_n) \rightarrow 1$.  

The following theorem whose proof is given in the Appendix shows that the robust estimators of the canonical directions   given in \eqref{eq:FCCest}  are consistent.

\begin{theorem}\label{consistency} 
	Let $(X_1,Y_1)\trasp,\dots, (X_n,Y_n)\trasp$ be i.i.d. with the same distribution as $(X,Y)\trasp\sim P_{(X,Y)}$. Assume that \ref{A1}-\ref{A2}, \ref{A3}(a) and \ref{A4} hold, then we have that
	\begin{enumerate}[font=\bfseries, label = (\alph*), ref=(\textbf{\alph*})]
		\item\label{itm:a}  $\wlam_{\tau}\convpp \lambda_0$, so the estimate of the canonical correlation is consistent,
		\item\label{itm:b} $\itL_{\rob}(\wphi_{\tau,1},\wpsi_{\tau,1})\convpp \lambda_0$,
		\item\label{itm:c} $\Lr^X(\wphi_{\tau,1},\phi_1) \convpp 1$  and $\Lr^Y(\wpsi_{\tau,1},\psi_1) \convpp 1$.
	\end{enumerate}
\end{theorem}

In order to get consistency results for the robust smoothed  canonical correlation estimators in the basis expansion domain, it is necessary to adopt additional notations and assumptions. Let $\bkapa=\bkapa_n=(\tau_n,d_n)$, $\wtlam_{\bkapa}  = \sup_{u,v\,\in\itH_d^0}\wLtr(u, v)$ $\lambda_{\bkapa}  = \sup_{u,v\,\in\itH_d^0}\Ltr(u, v)$, with solutions $(\wtphi_{\bkapa,1},\wtpsi_{\bkapa,1})$ and $(\phi_{\bkapa,1},\psi_{\bkapa,1})$, respectively. Let
\begin{align*}
D_{n,X}&=\sup_{u\in\itH_{d_n},\;\|{u}\|_{1,\tau{_n}}=1}\left|s_{n,X}^2(u)-\sigma^2_X(u)\right|\,, \\ 
D_{n,Y} &=\sup_{v\in\itH_{d_n},\;\|{v}\|_{2,\tau{_n}}=1}\left|s_{n,Y}^2(v)-\sigma^2_Y(v)\right|\, ,  \\
D_{n,XY} & =\sup_{u,v\in\itH_{d_n},\;\|{u}\|_{1,\tau{_n}}=\|{v}\|_{2,\tau{_n}}=1}\left|g_{n}(u,v)-\gamma_{XY}(u,v)\right|\,,
\end{align*}
and let $\Pi_{ \itH_{d_n}}$ be the orthogonal projection operator onto $ \itH_{d_n}$	
\begin{enumerate}[label = \textbf{C\arabic*}]
	\setcounter{enumi}{4}	
	\item\label{A5}   The basis  $\{\xi_i\}_{i\geq 1}\subset\HSo$ and $d=d_n$ is such that $d_n\rightarrow \infty$. The smoothing parameter $\tau=\tau_n\geq 0$ is such that  $\tau_n\rightarrow 0$,  $ \max(D_{n,X},D_{n,Y})\convpp 0$ and one of the following hold
	\begin{enumerate}
		\item  $ D_{n,XY}\convpp 0$ as $n\rightarrow \infty$ and there exists a constant $A>0$ such that for any $u,v\in \HSo$, we have that $$  {g_{n}^2(u,v)}/\left(\{s_{n,X}^2(u)+\tau\Psi(u)\}\{s_{n,Y}^2(v)+\tau\Psi(v)\}\right) \le A.$$
		\item The co--association measure is such that $g_{n}(u,v)=r_n(u,v) s_{n,X} (u)s_{n,Y}(v)$ and we have that  $\delta_n=\sup_{u,v\in\Hdo,\;\|{u}\|_{1,\tau{_n}}=\|{v}\|_{2,\tau{_n}}=1}\left|r_{n}(u,v)-\rho_{XY}(u,v)\right|$.
	\end{enumerate} 
	
	\item\label{A6} $\sigma^2_X(u):\itH\longrightarrow\real$, $\sigma^2_Y:\itH\longrightarrow\real$ and $\gamma_{XY}:\itH\times\itH\longrightarrow\real$ are continuous in $\phi_1$, $\psi_1$ and $(\phi_1,\psi_1)$, respectively.  
\end{enumerate}
Note that assumption \ref{A5}  is slightly weaker than \ref{A4}.

The following theorem whose proof is given in the Appendix shows that the robust estimators of the canonical directions  given in \eqref{eq:FCCestsieves}  are consistent. 

\begin{theorem}\label{consistencysieves} 
	Let $(X_1,Y_1)\trasp,\dots, (X_n,Y_n)\trasp$ be i.i.d. with the same distribution as $(X,Y)\trasp\sim P_{(X,Y)}$. Assume that \ref{A1}-\ref{A3}(a) and \ref{A5}-\ref{A6}  hold, $\tau_n\Psi(\Pi_{ \itH_{d_n}}\phi_1)\rightarrow 0$ and $\tau_n\Psi(\Pi_{ \itH_{d_n}}\psi_1)\rightarrow 0$, then we have that
	\begin{enumerate}[font=\bfseries, label = (\alph*), ref=(\textbf{\alph*})]
		\item\label{itm:as}  $\wtlam_{\bkapa}\convpp \lambda_0$, so the estimate of the canonical correlation is consistent,
		\item\label{itm:bs} $\itL_{\rob}(\wtphi_{\bkapa,1},\wtpsi_{\bkapa,1})\convpp \lambda_0$,
		\item\label{itm:cs} $\Lr^X(\wtphi_{\bkapa,1},\phi_1) \convpp 1$  and $\Lr^Y(\wtpsi_{\bkapa,1},\psi_1) \convpp 1$.
	\end{enumerate}
\end{theorem}

\subsection{Some general comments}
Assumptions \ref{A2} and \ref{A3} are similar to assumptions 3 and 4 in \cite{Leurgans}. In particular,  \ref{A3}(b) corresponds to the first part of assumption 4 in \cite{Leurgans}. \ref{A3} is satisfied, for example, when the roughness penalty is the integrated squared second derivative subject to periodic boundary conditions. Note that, the assumptions  $\tau\Psi(\Pi_{ \itH_{d_n}}\phi_1)\rightarrow 0$ and $\tau\Psi(\Pi_{ \itH_{d_n}}\psi_1)\rightarrow 0$ appearing in Theorem \ref{consistencysieves} are satisfied when $\{\Psi(\Pi_{ \itH_{d}}\phi_1)\}_{d\in\natu}$ and $\{\Psi(\Pi_{ \itH_{d}}\psi_1)\}_{d\in\natu}$ are bounded. 

Recall that a desirable property is that  the measures of co-association and scale defining $\itL_{\rob}$ determine the same canonical directions, which are the target ones, at least for a given distribution family. This property know as Fisher-consistency is strongly connected with \ref{A1} and \ref{A2}. In particular,  if  $\gamma_\rob$ is related to an association measure $\rho_\rob$  and the scale functional $\sigma_\rob$ by $\gamma_\rob \left( U,V \right)= \rho_\rob\left( U,V \right)\;\sigma_\rob \left( U\right)\sigma_\rob\left( V\right)$, then $ \itL_{\rob}(u,v)= \rho_{\rob}^2\left(P[ \langle u, X\rangle,\, \langle v, Y\rangle ]\right)$, so that \ref{A2} will be a consequence of the Fisher--consistency of  $\rho_{\rob}$ as discussed below. Some examples of association measures Fisher--consistent for elliptical distributed vectors were discussed in Remark \ref{remark:remark2}.

When $\gamma_\rob$ is the covariance and $\sigma_\rob$ is the standard deviation $\sd$, \ref{A1} holds, with $\bGa_{11}=\bGa_{XX}$, $\bGa_{22}=\bGa_{YY}$,  $\bGa_{12}=\bGa_{XY}$ and $c=1$. Notice that, in this case,  a necessary condition for a good definition of the canonical weights is that both random elements $X$ and $Y$  have finite second moments. This moment requirement may be relaxed when other association measures are considered.   

Assume that $(X,Y)\trasp\sim \itE(\mu,\bGa, \varphi)$ where $\itE(\bmu,\bGa, \varphi)$ denotes an elliptical distribution, as defined in \cite{Bali}, with parameters $\mu=(\mu_1,\mu_2)\trasp\in \itH\times \itH$  and  $\bGa$ is as in \eqref{eq:Gamma}.
As a consequence of the definition of elliptical random elements, for any $u,v\in\itH$, $(\langle u, X\rangle, \langle v, Y\rangle)\trasp \sim \itE_2(\bmu_{u,v},\bSi_{ u,v }, \varphi)$, where $\bmu_{u,v}= (\langle u,\mu_1\rangle, \langle v,\mu_2\rangle)\trasp$ and the  diagonal elements of $\bSi_{ u,v }$ are  $\langle u,\bGa_{11}u\rangle$ and $\langle v,\bGa_{22}v\rangle$, while the cross-diagonal ones are $\langle u,\bGa_{12}v\rangle$ and $\langle v,\bGa_{21}u\rangle\,$. As mentioned in Remark \ref{remark:remark2}, when considering a robust scale functional $\sigma_{\rob}(\cdot)$,  there exists a constant $c>0$ such that  $\sigma^2_X(u)=c\left\langle u,\bGa_{11} u\right\rangle$ and $\sigma^2_Y(v)=c\left\langle v,\bGa_{22} v\right\rangle$, for any $u,v\in \itH$, as required in  \ref{A1}. Let us denote $\rho_\rob$ the measure of association   $\rho_\rob \left( U,V \right)= \gamma_\rob\left( U,V \right)/\{\sigma_\rob \left( U\right)\sigma_\rob\left( V\right)\}$. If   $\rho_\rob$  is Fisher-consistent at the family of bivariate elliptical distributions,   we have that $\itL_{\rob}(u,v)=\rho_\rob^2(P_{(X,Y)}[u,v])=\itL_{\bGach}(u,v)$,  where
\begin{equation}
\itL_{\bGach}(u,v)=  \frac{\left\langle u, \bGa_{12}v\right\rangle^2}{\left\langle u, \bGa_{11}u\right\rangle\left\langle v, \bGa_{22}v\right\rangle}\,,
\label{eq:Lgama}
\end{equation}
so that $\gamma_{XY}(u,v)=c\left\langle u,\bGa_{12} v\right\rangle$ which corresponds to the  representation in \ref{A1}. 

A direct consequence of \eqref{eq:Lgama}, is that  the estimated canonical directions defined in \eqref{eq:FCCest} are consistent for the first canonical directions associated to $\bGa$, that is, to the maximizers of $\itL_{\bGach}(u,v)$,  which does not depend on $\rho_{\rob}$,  a fact that was already mentioned in \cite{Alvarez} and that states that the robust canonical directions are Fisher-consistent at  elliptical processes.
In particular,  if the process has second moments, taking into account that there exists an $a\in \real$, $a>0$,  such that the covariance operator of $(X,Y)\trasp$ equals ${a}\bGa$, we have that $\itL_{\rob}(u,v)=\itL_{\bGach}(u,v)=\itL(u,v) $, so the robust functionals related to the canonical analysis are the usual ones.
As  mentioned in Remark \ref{remark:remark2}, all association measures in Section \ref{measures} are Fisher-consistent for bivariate elliptical families. Hence, the above discussion implies  that    \ref{A2} holds   for  these association measures if the conditions in Theorem 4.8 of \cite{He} hold. 

The following Lemma gives conditions ensuring that the convergences in \ref{A4} hold, i.e., that $ C_{n,X}\convpp 0$, $C_{n,Y}\convpp 0$ and  $ C_{n,XY}\convpp 0$.

\begin{lemma}\label{lema:CnX} 
	Let $(X_1,Y_1)\trasp,\dots, (X_n,Y_n)\trasp$ be i.i.d. with the same distribution as $(X,Y)\trasp\sim P_{(X,Y)}$. Let $\eta_n\;=\;\sup_{\|{u}\|=\|{v}\|=1}\left|g_{n}(u,v)-\gamma_{XY}(u,v)\right|\,$,\hspace{0.8cm} $\zeta_n=\sup_{\|{u}\|=1}\left|s_{n,X}^2(u)-\sigma^2_X(u)\right|$ and  
	$\nu_n=\sup_{\|{v}\|=1}\left|s_{n,Y}^2(v)-\sigma^2_Y(v)\right|$.  If  \ref{A1} and \ref{A3} hold,  $\tau_n\rightarrow 0$ and    $\tau_n^{-1}\max(\zeta_n,\nu_n,\eta_n)\convpp 0$ as $n\rightarrow \infty$, then we have that $ C_{n,X}\convpp 0$, $C_{n,Y}\convpp 0$ and  $ C_{n,XY}\convpp 0$.
\end{lemma}

It is worth noticing that Lemma \ref{lema:CnX} and Theorem \ref{consistency} allow to derive strong consistency results for the canonical directions defined in \cite{Leurgans}. 
Effectively, define $L(t)=\log\;\max(t, e)$ and $LL(t)=L(L(t))$, for any $t>0$. Moreover, we will denote $LLn= LL(n)$, so that $LLn=\log\:\log n$ for $n\ge 3$ and $LLn=1$ for $n=1,2$. 
Let  $Z=(X,Y)\trasp $  and  $\bGa_{ZZ}=\esp\left[\{Z-\esp(Z)\}\otimes \{Z-\esp(Z)\}\right]$ its covariance operator. Note that  $\bGa_{ZZ}$ is a self--adjoint continuous linear operator over $\itH\times \itH$; moreover, it is a Hilbert-Schmidt operator. For simplicity, $\itF$ will stand for the Hilbert space of such operators with inner product defined by  $\langle\bGa_1, \bGa_2\rangle_{\itF} = \mbox{trace}(\bGa_1^{{*}}\bGa_2) =\sum_{j=1}^\infty \langle\bGa_1 z_j, \bGa_2 z_j\rangle_{{\itH\times \itH}}$, where $\{z_j : j\ge 1\}$ is any orthonormal basis of $\itH\times \itH$ {and $\bGa_1^*$ is the adjoint of $\bGa_1$}. 
Furthermore, define $V=\{Z-\esp(Z)\}\otimes \{Z-\esp(Z)\}- \bGa_{ZZ}$ which is a zero mean random  element in $\itF$. Then, if $\esp\left\{\|V\|_{\itF}^{ 2}/LL(\|V\|_{\itF})\right\}<\infty$ and $\esp\left(\langle V, F\rangle_{\itF}^{ 2}\right)<\infty$, for any $F\in \itF$, the law of iterated logarithm in Hilbert spaces obtained in \cite{Acosta} allows us to conclude  that the assumptions in Lemma \ref{lema:CnX} hold when $\tau_n\; \sqrt{n/LLn}\to \infty$. Hence, under \ref{A2} {and \ref{A3}}, the canonical directions are consistent in the $\bGa_{ZZ}-$norm.

\section{Conclusion}  \label{sec:conclusion}

We have introduced two procedures to obtain  robust estimators of the canonical directions based on co--association measures and a regularization term involving a roughness penalty. The resulting estimators are consistent under mild conditions. Furthermore, if the process $(X,Y)\trasp$ has an elliptical distribution with finite second moment, the resulting target quantities correspond to the usual canonical correlation and directions.   We have assumed that $X$ and $Y$ are defined over the same infinite--dimensional space $\itH$ to avoid burden notation. The extension to the situation in which $X\in \itH_1$ and $Y\in \itH_2$ is straightforward.

\vskip 14pt
\noindent {\bf Acknowledgements.}
{\small   This research was partially supported by Grants \textsc{pict} 2018-00740 from \textsc{anpcyt},   20020170100022BA    from the Universidad de Buenos Aires  at Buenos Aires, Argentina (Graciela Boente), the Spanish Project {MTM2016-76969P} from the Ministry of Economy	and Competitiveness, Spain (MINECO/ AEI/FEDER, UE) (Graciela Boente) and   \textsc{ppid} \textsc{x}030 and \textsc{pid} \textsc{i}231 from Universidad Nacional de La Plata, Argentina  (Nadia Kudraszow).} 


\setcounter{section}{0}
\renewcommand{\thesection}{A}

\setcounter{equation}{0}
\renewcommand{\theequation}{A.\arabic{equation}}

\section{Appendix}{\label{sec:appen}}
To prove Theorem \ref{consistency}, we need some preliminary results. 
\begin{lemma}\label{lemma3}
	Assume that \ref{A2}  holds. Then,  for any $u$ and $v$ in $\itH_{\smooth}$, we have that $\Lr(u, v) \geq \Ltr(u, v)$ and $\lambda_{\tau}\leq\Lr(\phi_{\tau,1},\psi_{\tau,1})$. Moreover,  $\lambda_{\tau}\rightarrow\lambda_0$ as $\tau\rightarrow 0$.
\end{lemma}

\noindent\textbf{Proof.} To show that $\Lr(u, v) \geq \Ltr(u, v)$ note that, given any $u$ and $v$ in $\itH_{\smooth}$
\begin{equation}\label{LsobreL}
\frac{\Ltr(u,v)}{\Lr(u,v)}=\frac{\sigma^2_X(u)}{\{\sigma^2_X(u)+\tau\Psi(u)\}}\frac{\sigma^2_Y(v)}{\{\sigma^2_Y(v)+\tau\Psi(v)\}}\leq 1
\end{equation}
where the inequality follows from the fact that $\sigma^2_\rob$, $\tau$ and $\Psi(\cdot)$ are non negative. Thus, $\lambda_{\tau}=\Ltr(\phi_{\tau,1},\psi_{\tau,1})\leq\itL_{\rob}(\phi_{\tau,1},\psi_{\tau,1})\leq\lambda_0$, proving the second part. From the expression for the ratio $\Ltr(u,v)/\Lr(u,v)$ given  in (\ref{LsobreL}) we conclude  that, for any fixed $u$ and $v$ in $\itH_{\smooth}$, $\Ltr(u, v)\rightarrow\Lr(u, v)$ as $\tau\rightarrow0$. Hence, using that $\lambda_0 \geq \lambda_{\tau}\geq \Ltr(\phi_{1},\psi_{1})$ and the fact that $\Ltr(\phi_{1},\psi_{1})\rightarrow\Lr(\phi_{1},\psi_{1}) = \lambda_0$, we conclude the proof. \qed

The following result  provides a sufficient condition for convergence in the $\Lr$-norm. It may be derived using similar arguments to those considered in  Lemma 4 in \cite{Leurgans}, for that reason, its proof is omitted.

\begin{lemma}\label{lemma4}
	Assume that \ref{A1} and \ref{A2} hold and that $\Lr(u_n, v_n) \rightarrow \rho_0^2=\itL_{\rob}(\phi_{1},\psi_{1})$ as $n \rightarrow \infty$. Then $(u_n, v_n) \rightarrow (\phi_{1},\psi_{1})$ in the $\Lr$-norm.
\end{lemma}

The next proposition is the key step in the  proof of Theorem \ref{consistency}.

\begin{proposition}\label{prop3}
	Assume that  \ref{A3}(a) and \ref{A4}  hold, then $\sup_{u,v\,\in\HSo} |\wLtrn(u,v)-\Ltrn(u,v)|\convpp 0 $. 
\end{proposition}

\noindent\textsc{Proof.} From the equivariance of the scale functional and \ref{A4}, it follows that
$$\sup_{u\,\in\HSo}\left|\frac{s_{n,X}^2(u)+\tau\Psi(u)}{\sigma^2_X(u)+\tau\Psi(u)}-1\right|=C_{n,X}\convpp 0\,.$$
Similarly, we have that
$$
\sup_{v\,\in\HSo}\left|\frac{s_{n,Y}^2(v)+\tau\Psi(v)}{\sigma^2_Y(v)+\tau\Psi(v)}-1\right|=C_{n,Y}\convpp 0\,.
$$
Note that, if \ref{A4}(a) holds,  
\begin{eqnarray*}
\sup_{u,v\,\in\HSo}\frac{\left|g_{n}(u,v)-\gamma_{XY}(u,v)\right|}{\left\{\sigma^2_{X}(u) \rangle+\tau\Psi(u)\right\}^{1/2}\, 
	\left\{\sigma^2_{Y}(v)+\tau\Psi(v)\right\}^{1/2}}
&=&\sup_{u, v\,\in\HSo}\frac{\left|g_{n}(u,v)-\gamma_{XY}(u,v)\right| }{\|{u}\|_{1,\tau} \|{v}\|_{2,\tau} }
\end{eqnarray*}
\begin{eqnarray*}
\hspace{18em}&=&   \sup_{\|{u}\|_{1,\tau}=\|{v}\|_{2,\tau}=1}\left|g_{n}(u,v)-\gamma_{XY}(u,v)\right|\\
&=&C_{n,XY}\convpp 0\,.
\end{eqnarray*}
Let us first prove our assertion when    \ref{A4}(a) holds. Note that
\begin{align*}
&|\wLtrn(u,v)-\Ltrn(u,v)| = \\
&= \left| \frac{g_{n}^2(u,v) }{\{s_{n,X}^2(u)+\tau\Psi(u)\}\{s_{n,Y}^2(v)+\tau\Psi(v)\}}- \frac{\gamma_{XY}^2(u,v) }{\{\sigma^2_{X}(u)+\tau\Psi(u)\}\{\sigma^2_{Y}(v)+\tau\Psi(v)\}}\right|\\
&\le  \frac{g_{n}^2(u,v)}{\{s_{n,X}^2(u)+\tau\Psi(u)\}\{s_{n,Y}^2(v)+\tau\Psi(v)\}} \left| \frac{\{s_{n,X}^2(u)+\tau\Psi(u)\}\{s_{n,Y}^2(v)+\tau\Psi(v)\}}{\{\sigma^2_{X}(u)+\tau\Psi(u)\}\{\sigma^2_{Y}(v)+\tau\Psi(v)\}}- 1\right|\\
&+ \frac{|g_{n}^2(u,v)-\gamma_{XY}^2(u,v) |}{\{\sigma^2_{X}(u)+\tau\Psi(u)\}\{\sigma^2_{Y}(v)+\tau\Psi(v)\}} =A_{1,n}(u,v)+A_{2,n}(u,v)
\end{align*}
Let us begin by bounding $A_{2,n}(u,v)$. Using that for any $u, v\in \HSo$, we have that 
$$\frac{ |\gamma_{XY}^2(u,v) |}{ \{\sigma^2_{X}(u)+\tau\Psi(u)\}\{\sigma^2_{Y}(v)+\tau\Psi(v)\} }\le\itL_{\rob}(u, v)\le  \lambda_0\,,$$ 
and using   \ref{A4}(a), we obtain that
\begin{align*}
\sup_{u,v\,\in\HSo} A_{2,n} (u,v)&\le \sup_{u,v\,\in\HSo} \frac{\left\{g_{n} (u,v)-\gamma_{XY} (u,v) \right\}^2 + 2 |\gamma_{XY} (u,v)|\; \left|g_{n} (u,v)-\gamma_{XY} (u,v) \right|}{\{\sigma^2_{X}(u)+\tau\Psi(u)\}\{\sigma^2_{Y}(v)+\tau\Psi(v)\}}\\
&\le C_{n,XY}^2+ \lambda_0^{1/2} \; C_{n,XY}\convpp 0\,.
\end{align*}
On the other hand, to bound $A_{1,n}(u,v)$, we use that $a  b  -1= (a  -1)(b-1)+(a-1)+(b-1)$. Choosing $a= \{s_{n,X}^2(u)+\tau\Psi(u)\}/\{\sigma^2_{X}(u)+\tau\Psi(u)\}$ and $b=\{s_{n,Y}^2(v)+\tau\Psi(v)\}/ \{\sigma^2_{Y}(v)+\tau\Psi(v)\}$, we obtain that
\begin{eqnarray*}
\sup_{u,v\,\in\HSo}  &&\left| \frac{\{s_{n,X}^2(u)+\tau\Psi(u)\}\{s_{n,Y}^2(v)+\tau\Psi(v)\}}{\{\sigma^2_{X}(u)+\tau\Psi(u)\}\{\sigma^2_{Y}(v)+\tau\Psi(v)\}}-1\right| \le\\
&&\hspace{4.5em}\le\sup_{u,v\,\in\HSo}  \left|\frac{ s_{n,X}^2(u)+\tau\Psi(u) }{ \sigma^2_{X}(u)+\tau\Psi(u) }- 1  \right|\; \sup_{u,v\,\in\HSo}\left|\frac{ s_{n,Y}^2(v)+\tau\Psi(v) }{ \sigma^2_{Y}(v)+\tau\Psi(v) }- 1\right|\\
&&\hspace{4.5em}+ \sup_{u,v\,\in\HSo}  \left|\frac{ s_{n,X}^2(u)+\tau\Psi(u) }{ \sigma^2_{X}(u)+\tau\Psi(u) }- 1  \right|+ \sup_{u,v\,\in\HSo}\left|\frac{ s_{n,Y}^2(v)+\tau\Psi(v) }{ \sigma^2_{Y}(v)+\tau\Psi(v) }- 1\right|\\
&&\hspace{4.5em} \le C_{n,X}\; C_{n,Y}+  C_{n,X}+ C_{n,Y} \convpp 0\;.
\end{eqnarray*}
Using that, for any $u,v\,\in\HSo$,
$$ \displaystyle{\frac{g_{n}^2(u,v)}{\{s_{n,X}^2(u)+\tau\Psi(u)\}\{s_{n,Y}^2(v)+\tau\Psi(v)\}} }\le A$$
we get that $\sup_{u,v\,\in\HSo} A_{1,n}(u,v)\convpp 0$, which concludes the proof when  \ref{A4}(a) holds.

Assume now that   \ref{A4}(b) holds. Then, $\gamma_n^2(u,v)\;=\;r_{n}^2(u,v)s_{n,X}^2(u)s_{n,Y}^2(v)\;$, $\gamma_{XY}^2(u,v)\;=\;\rho_{XY}^2(u,v)\sigma^2_{X}(u)\sigma^2_{Y}(v)\;$. Thus,  for any
 ${u,v\,\in\HSo}$, $$|\wLtrn(u,v)-\Ltrn(u,v)| \le B_{1,n}(u,v)+B_{2,n}(u,v)$$ where
\begin{align*}
&B_{1,n}(u,v)
 =  \frac{\left| r_{n}^2(u,v)-\rho_{XY}^2(u,v)\right|s_{n,X}^2(u)s_{n,Y}^2(v)}{\{s_{n,X}^2(u)+\tau\Psi(u)\}\{s_{n,Y}^2(v)+\tau\Psi(v)\}}\qquad \;\;\;\;\text{ and }  \;\;\;\;\qquad
B_{2,n}(u,v)  = \\ &\left|\frac{s_{n,X}^2(u)s_{n,Y}^2(v)}{\{s_{n,X}^2(u)+\tau\Psi(u)\}\{s_{n,Y}^2(v)+\tau\Psi(v)\}}- \frac{\sigma^2_{X}(u)\sigma^2_{Y}(v)}{\{\sigma^2_{X}(u)+\tau\Psi(u)\}\{\sigma^2_{Y}(v)+\tau\Psi(v)\}}\right|.
\end{align*}
Using that  $r_{n}^2(u,v)\le 1$ and  $\rho_{XY}^2(u,v)\le 1$, from \ref{A4}(b) we immediately obtain that 
$B_{1,n}\le   2\; \sup_{\|{u}\|=\|{v}\|=1}\left|r_{n}(u,v)-\rho_{XY}(u,v)\right|= 2\; \theta_n\convpp 0\,.$
It only remains to show that $\sup_{u,v\,\in\HSo}B_{2,n} (u,v)\convpp 0$. Note that $\sup_{u,v\,\in\HSo}B_{2,n} (u,v)$ can be bounded by
\begin{align*}
\sup_{u,v\,\in\HSo}B_{2,n} (u,v)
& \le  C_{n,Y}  +\sup_{u,v\,\in\HSo}\left|\frac{s_{n,X}^2(u)s_{n,Y}^2(v)}{\{s_{n,X}^2(u)+\tau\Psi(u)\}\{\sigma^2_{Y}(v)+\tau\Psi(v)\}}-\right.\\ &\hspace{10em}-\left.\frac{\sigma^2_{X}(u)\sigma^2_{Y}(v)}{\{\sigma^2_{X}(u)+\tau\Psi(u)\}\{\sigma^2_{Y}(v)+\tau\Psi(v)\}}\right|\,.
\end{align*}
Denote $D_{1,n}$ the second term on the  right hand side of the above equation. Arguing similarly we get that
\begin{align*}
D_{1,n} 
& \le   C_{n,X}\; \sup_{u,v\,\in\HSo} \frac{s_{n,Y}^2(v)}{ \sigma^2_{Y}(v)+\tau\Psi(v)}  +\sup_{u,v\,\in\HSo}\frac{\left|s_{n,X}^2(u)s_{n,Y}^2(v)-\sigma^2_{X}(u)\sigma^2_{Y}(v)\right|}{\{\sigma^2_{X}(u)+\tau\Psi(u)\}\{\sigma^2_{Y}(v)+\tau\Psi(v)\}}\\
& \le   C_{n,X}+ C_{n,X}\;\sup_{u,v\,\in\HSo} \frac{|s_{n,Y}^2(v)-\sigma^2_{Y}(v)|}{ \sigma^2_{Y}(v)+\tau\Psi(v)}+\\
& +\sup_{u,v\,\in\HSo}\frac{\left|s_{n,X}^2(u)s_{n,Y}^2(v)-\sigma^2_{X}(u)\sigma^2_{Y}(v)\right|}{\{\sigma^2_{X}(u)+\tau\Psi(u)\}\{\sigma^2_{Y}(v)+\tau\Psi(v)\}}\\ 
& \le C_{n,X}+ C_{n,X}\;C_{n,Y} +\sup_{u,v\,\in\HSo}\frac{\left|s_{n,X}^2(u)s_{n,Y}^2(v)-\sigma^2_{X}(u)\sigma^2_{Y}(v)\right|}{\{\sigma^2_{X}(u)+\tau\Psi(u)\}\{\sigma^2_{Y}(v)+\tau\Psi(v)\}}
\\
&=  C_{n,X}+ C_{n,X}\;C_{n,Y}+D_{2,n} \,.
\end{align*}
Finally $D_{2,n} $ can easily be bounded by $C_{n,X}\, C_{n,Y}+C_{n,X}+C_{n,Y}$ so $$
\sup_{u,v\,\in\HSo}B_{2,n} (u,v) \le C_{n,Y}  + 2\left(C_{n,X}+ C_{n,X}\;C_{n,Y}\right)\convpp 0\;.
$$
concluding the proof.  \qed

\bigskip

\noindent\textsc{Proof of Theorem \ref{consistency}.} (a) Standard arguments  and Proposition \ref{prop3} allow to conclude that $|\wlam_{\tau}-\lambda_{\tau}|\convpp 0$ which together with Lemma \ref{lemma3} entails that $\wlam_{\tau}\convpp\lambda_0$. 

\noi To prove (b) first note  that, by Proposition \ref{prop3}, $
|\wLtr(\wphi_{\tau,1},\wpsi_{\tau,1})-\Ltr(\wphi_{\tau,1},\wpsi_{\tau,1})|\convpp 0$.
Using $\sim_{a.s.}$ to connect quantities whose difference converges to 0 almost surely, from Lemma \ref{lemma3}, we get  that 
$$
\lambda_0\geq\Lr(\wphi_{\tau,1},\wpsi_{\tau,1})\geq\Ltr(\wphi_{\tau,1},\wpsi_{\tau,1})\sim_{a.s.}\wLtr(\wphi_{\tau,1},\wpsi_{\tau,1})=\wlam_{\tau}\convpp \lambda_0\,,
$$
concluding the proof of (b). 

\noi (c) is  an immediate consequence of (b) applying Lemma \ref{lemma4}.  \qed

\bigskip

\noindent\textsc{Proof of Theorem \ref{consistencysieves}.} (a) Using \ref{A3}(a) and \ref{A5}, the following result may be derived using similar arguments to those considered in the proof of  Proposition \ref{prop3},
\begin{equation}\label{P1sieves}
\sup_{u,v\,\in\Hdo} |\wLtrn(u,v)-\Ltrn(u,v)|\convpp 0\;,
\end{equation}
which implies  $|\wtlam_{\bkapa}-\lambda_{\bkapa}|\convpp 0$. From Lemma \ref{lemma3}, it is easily seen that  $\lambda_0 \geq \lambda_{\tau}\geq \lambda_{\bkapa}\geq \Ltr\left(\wtphi_{d_n,1},\wtpsi_{d_n,1}\right)$, where $\wtphi_{d_n,1}=  {\Pi_{ \itH_{d_n}}\phi_1}/{\|\Pi_{ \itH_{d_n,1}} \phi_1\|}$ and $\wtpsi_{d_n,1}={\Pi_{ \itH_{d_n}}\psi_1}/{\|\Pi_{ \itH_{d_n}} \psi_1\|}$ are the standardized orthogonal projections of $\phi_1$ and $\psi_1$ respectively, onto $ \itH_{d_n}$. Then, the proof of (a) is completed by showing that $$\Ltr\left(\wtphi_{d_n,1},\wtpsi_{d_n,1}\right)\rightarrow \Lr(\phi_{1},\psi_{1})=\lambda_0.$$ Since  $\|\Pi_{ \itH_{d_n}} \phi_1\|\rightarrow\| \phi_1\|$,  $\|\Pi_{ \itH_{d_n}} \psi_1\|\rightarrow\| \psi_1\|$, $\tau\Psi(\Pi_{ \itH_{d_n}}\phi_1)\rightarrow 0$ and $\tau\Psi(\Pi_{ \itH_{d_n}}\psi_1)\rightarrow 0$, by \ref{A6},  we have
\begin{align*}
&\frac{\Ltr\left(\wtphi_{d_n,1},\wtpsi_{d_n,1}\right)}{\Lr(\phi_{1},\psi_{1})}=\\
&\hspace{3em}=\frac{\gamma_{XY}^2\left(\wtphi_{d_n,1},\wtpsi_{d_n,1}\right)}{\gamma_{XY}^2(\phi_{1},\psi_{1})}\frac{\sigma^2_X(\phi_{1})}{\left\{\sigma^2_X\left(\wtphi_{d_n,1}\right)+\tau\Psi\left(\wtphi_{d_n,1}\right)\right\}}\frac{\sigma^2_Y(\psi_{1})}{\left\{\sigma^2_Y\left(\wtpsi_{d_n,1}\right)+\tau\Psi\left(\wtpsi_{d_n,1}\right)\right\}}
\end{align*}
tends to 1.

Finally, (b) and (c) may be obtained in a similar fashion as (b) and (c) of Theorem \ref{consistency}.\qed 

\bigskip
Lemma \ref{lemma2}  will be needed in the proof of Lemma \ref{lema:CnX}. It corresponds to Lemma 2 in \cite{Leurgans} so its proof is omitted. 

\begin{lemma}\label{lemma2} 
	Assume that \ref{A1}, \ref{A3}(a)  and \ref{A3}(c) hold and let $\ell_1(\tau)$ and $\ell_2(\tau)$ be the smallest eigenvalues of $ c\bGa_{11}+\tau\Psi(\cdot) $ and $ c\;\bGa_{22}+\tau\Psi(\cdot)$,  respectively. Then, for $0<\tau \leq 1$, we have that $\ell_1(\tau)\geq \tau \ell_1(1)>0$ and $\ell_2(\tau)\geq \tau \ell_2(1)>0$.
\end{lemma}
\vskip0.1in
\noi \textsc{Proof of Lemma \ref{lema:CnX}}. Note that
\begin{align*}
C_{n,X} &= \sup_{\|{u}\|_{1,\tau_n}=1}\left|s_{n,X}^2(u)-\sigma^2_X(u)\right|
=\sup_{u\,\in\HSo}\left|\frac{s_{n,X}^2(u)-\sigma^2_X(u)}{\|{u}\|_{1,\tau_n}^2}\right| \\
& = \sup_{u\,\in\HSo}\left|\frac{s_{n,X}^2(u)-\sigma^2_X(u)}{\left\langle u,\bGa_{11}u\right\rangle+\tau_n\Psi(u)}\right| 
\leq \frac{\sup_{ u\,\in\HSo, \|u\|=1}\left|s_{n,X}^2(u)-\sigma^2_X(u)\right|}{\inf_{u\,\in\HSo, \|u\|=1}\{\left\langle u, \bGa_{11}u\right\rangle+\tau_n\Psi(u)\}}\;.
\end{align*}
Besides, if as in Lemma \ref{lemma2}, $\ell_1(\tau)$ stands for the smallest eigenvalue  of $  \bGa_{11}+\tau\Psi(\cdot) $, we conclude that for any $u \in\HSo$, $ \left\langle u, \bGa_{11}u\right\rangle+\tau\Psi(u) \ge \ell_1(\tau) \|u\|^2$, so that
$\inf_{u\,\in\HSo, \|u\|=1}\{ \left\langle u,\bGa_{11}u\right\rangle+\tau\Psi(u)\} \ge \ell_1(\tau)\,.$ 
Hence,   we get that
\begin{equation*}
C_{n,X}  \leq \frac{\sup_{u\,\in\HSo\,,\, \|u\|=1}\left|s_{n,X}^2(u)-\sigma^2_X(u)\right|}{\ell_1(\tau_n)\; }  \leq \frac{\zeta_n}{\ell_1(\tau_n)\; }\leq \frac{\zeta_n}{\tau_n \ell_1(1)}\convpp 0\,.
\end{equation*}
Similar arguments and the fact that  $\nu_n/\tau_n\convpp 0$, allow to show that $C_{n,Y} \convpp 0$. Finally,
\begin{align*}
C_{n,XY} & = \sup_{\|{u}\|_{1,\tau_n}=\|{v}\|_{2,\tau_n}=1}\left|g_{n}(u,v)-\gamma_{XY}(u,v)\right|
=\sup_{u, v\,\in\HSo}\left|\frac{g_{n}(u,v)-\gamma_{XY}(u,v)}{\|{u}\|_{1,\tau_n} \|{v}\|_{2,\tau_n} }\right| \\
&= \sup_{u,v\,\in\HSo}\left|\frac{g_{n}(u,v)-\gamma_{XY}(u,v)}{\left(\langle u,\bGa_{11}u \rangle+\tau_n\Psi(u)\right)^{1/2}\, 
	\left(\langle v,\bGa_{22}v \rangle+\tau_n\Psi(v)\right)^{1/2}}\right| \\
&\leq \frac{\sup_{u, v\,\in\HSo, \|u\|=\|v\|=1}\left|g_{n}(u,v)-\gamma_{XY}(u,v)\right|}
{\inf_{u\,\in\HSo, \|u\|=1}\{\left\langle u, \bGa_{11}u\right\rangle+\tau_n\Psi(u)\}^{1/2}\; \inf_{v\,\in\HSo, \|v\|=1}\{\left\langle v, \bGa_{22}v\right\rangle+\tau_n\Psi(v)\}^{1/2}}\,.
\end{align*}
Denoting as above,   $\ell_1(\tau)$  and $\ell_2(\tau)$ the smallest eigenvalues  of $  \bGa_{11}+\tau\Psi(\cdot) $ and    $  \bGa_{22}+\tau\Psi(\cdot) $, respectively, and using that from  Lemma \ref{lemma2} $\ell_j(\tau)\ge \tau \ell_j(1)$, for $j=1,2$, we obtain that 
\begin{align*}
C_{n,XY} & \leq \frac{\sup_{u, v\,\in\HSo, \|u\|=\|v\|=1}\left|g_{n}(u,v)-\gamma_{XY}(u,v)\right|}
{\left(\ell_1(\tau_n) \; \ell_2 (\tau_n) \right)^{1/2}}\le \frac{\eta_n}{\tau_n \;\left(\ell_1(1) \; \ell_2 (1) \right)^{1/2}}\,,
\end{align*}
which concludes the proof since $\eta_n/\tau_n\convpp 0$. \qed


\end{document}